\DeclareUnicodeCharacter{2500}{}
\documentclass[11pt]{article}
\usepackage{xurl} 
\usepackage{cite} 

\usepackage[utf8]{inputenc}
\usepackage[T1]{fontenc}
\usepackage{geometry}
\usepackage{lmodern,microtype}
\usepackage{amsmath,amssymb,amsfonts,amsthm,mathtools}
\usepackage{bm,bbm}
\usepackage{graphicx}
\usepackage{booktabs}
\usepackage{enumitem}
\usepackage{algorithm}
\usepackage{algpseudocode}
\usepackage{hyperref}
\usepackage[nameinlink,capitalize]{cleveref}
\usepackage{tikz-cd}
\usepackage{subcaption} 
\usepackage[font=small,labelfont=bf]{caption}
\usepackage[sort&compress,numbers]{natbib}
\usepackage{xcolor}   
\usepackage{lineno}   
\usepackage{amsmath, amssymb, amsthm, mathrsfs}

\title{\textbf{Persistent Sheaf Laplacian Analysis of Protein Stability and Solubility Changes upon Mutation}}
\author{Yiming Ren$^{1}$,  Junjie Wee$^{1}$, Xi Chen$^2$, Grace Qian$^3$, and Guo-Wei Wei$^{1,4,5}$\footnote{
			Corresponding author.		Email: weig@msu.edu} \\
		\\
		$^1$Department of Mathematics, \\
		Michigan State University, East Lansing, MI 48824, USA.\\
        $^2$The Frazer School, 4700 NW 89 Blvd, Gainesville, FL 32606, USA\\
        $^3$Lassiter High School, Marietta, GA 30066, USA\\
        $^4$Department of Biochemistry and Molecular Biology,\\
		Michigan State University, East Lansing, MI 48824, USA.  \\
		$^5$Department of Electrical and Computer Engineering,\\
		Michigan State University, East Lansing, MI 48824, USA. \\
        \\
	}

\begin{document}
\maketitle
\begin{abstract}

Genetic mutations frequently disrupt protein structure, stability, and solubility, acting as primary drivers for a wide spectrum of diseases. Despite the critical importance of these molecular alterations, existing computational models often lack interpretability, and fail to integrate essential physicochemical interaction. To overcome these limitations, we propose SheafLapNet, a unified predictive framework grounded in the mathematical theory of Topological Deep Learning (TDL) and Persistent Sheaf Laplacian (PSL). Unlike standard Topological Data Analysis (TDA) tools such as persistent homology, which are often insensitive to heterogeneous information, PSL explicitly encodes specific physical and chemical information such as partial charges directly into the topological analysis. SheafLapNet synergizes these sheaf-theoretic invariants with advanced protein transformer features and auxiliary physical descriptors to capture intrinsic molecular interactions in a multiscale and mechanistic manner. To validate our framework, we employ rigorous benchmarks for both regression and classification tasks. For stability prediction, we utilize the comprehensive S2648 and S350 datasets. For solubility prediction, we employ the PON-Sol2 dataset, which provides annotations for increased, decreased, or neutral solubility changes. By integrating these multi-perspective features, SheafLapNet achieves state-of-the-art performance across these diverse benchmarks, demonstrating that sheaf-theoretic modeling significantly enhances both interpretability and generalizability in predicting mutation-induced structural and functional changes.
\end{abstract}


\newpage

\section{Introduction}

Genetic mutations induce alterations in the amino acid sequence that can severely disrupt the delicate thermodynamic balance of protein structures. These atomic-level perturbations frequently compromise critical physicochemical attributes, including folding stability, binding affinity, reactivity, and solubility, often precipitating loss-of-function phenotypes or toxic aggregation linked to severe pathologies. Such molecular mechanisms are central to the etiology of neurodegenerative disorders like Alzheimer's \cite{hurley2023familial} and Parkinson's disease \cite{funayama2023molecular}, as well as various cancers \cite{chen2022mutant} and metabolic syndromes \cite{zhang2013y328c, zhang2011silico}. Despite the clinical criticality of these stability and solubility profiles, accurately predicting mutation-induced changes remains a formidable challenge. Protein solubility is governed by a multifaceted interplay of factors, ranging from intrinsic sequence motifs and post-translational modifications to extrinsic environmental conditions such as solvent type, ionic strength, and temperature. This complexity is compounded by the sheer volume of genomic data; while more than four million missense variants have been cataloged, the majority remain classified as variants of uncertain significance (VUS) \cite{karczewski2020mutational}, and existing datasets often lack the comprehensive environmental annotations required to resolve these nuances. Given that traditional experimental validation remains resource-intensive and low-throughput \cite{marian2020clinical, molotkov2024making}, there is an urgent need for advanced computational approaches capable of decoding these complex biophysical interactions to reliably predict the functional impact of genetic variants.

Protein stability and solubility are fundamental to protein function and directly involved in human disease, yet predicting mutation-induced changes in these properties remains a complex computational challenge. In the domain of stability, deep learning architectures have capitalized on large-scale datasets. For instance, mutDDG-SSM \cite{li2024prediction} hybridizes graph attention networks with gradient boosting trees, while TopologyNet \cite{cang2017topologynet} employs a multi-task framework to link stability data with disease associations. Despite these advances, models developed over the past decade, such as STRUM \cite{quan2016strum}, have leveraged modern machine learning to uncover hidden relationships between structure and stability but often provide limited interpretability. Parallel efforts in the prediction of mutation-induced solubility changes have yielded tools such as CamSol \cite{sormanni2015camsol}, SODA \cite{paladin2017soda}, Solubis \cite{van2016solubis}, PON-Sol \cite{yang2016pon}, and PON-Sol2\cite{yang2021pon}, which utilizes gradient boosting on expanded datasets. However, a unified framework capable of simultaneously addressing these coupled properties remains absent. Furthermore, despite their predictive power, many existing models may fail to explicitly account for fundamental physical interactions, including hydrogen bonding, van der Waals forces, hydrophobicity, and electrostatics, that govern molecular behavior \cite{sun2022electrostatics, stefl2013molecular}. This limitation, coupled with suboptimal performance metrics like the normalized Correct Prediction Ratio (CPR) in solubility tasks, underscores an urgent demand for Explainable AI (XAI) frameworks capable of providing interpretable, mechanistically grounded insights into disease causality.

To address these challenges, we turn to Topological Data Analysis (TDA), an emerging mathematical field that utilizes algebraic topology to analyze structural patterns within complex data. Its central tool, persistent homology \cite{zomorodian2004computing, epstein2011topological, carlsson2009topology}, integrates classical homology and filtration to create a multiscale analysis of data. TDA becomes a powerful approach in data science when it is paired with machine learning. In particular, topological Deep Learning (TDL) that integrates TDA and deep neural network was introduced for the first time in 2017 \cite{cang2017topologynet}. Over years, TDL has become a new frontier in rational learning  \cite{papamarkou2024position}. various deep neural networks have been proposed  \cite{barbero2022sheaf,hajij2025copresheaf}. 
TDL has achieved tremendous success in deciphering biomolecules \cite{cang2018representability}. 
Some of the most compelling application examples in which TDL has consistently demonstrated clear advantages over competing methods include its strong performance in the D3R Grand Challenges \cite{nguyen2019mathematical}, its role in uncovering the evolutionary mechanisms of SARS-CoV-2 \cite{chen2020mutations}, and its successful forecasting of viral  variants \cite{chen2022omicron}.

However, traditional TDA with persistent homology faces inherent limitations as it is insensitive to homotopic shape evolution without topological changes, cannot distinguish between geometric isomers, and is incapable of differentiating atom types and  encoding directed relations. To overcome these limitations, the persistent spectral theory \cite{wang2020persistent}, also known as the Persistent Laplacian (PL) \cite{memoli2022persistent, wang2021hermes}, was proposed. By capturing homotopic shape evolution via nonharmonic spectra, computational algorithms such as the HERMES software \cite{wang2021hermes}, homotopy continuation \cite{wei2021homotopy} and PETLS \cite{jones2025petls} have facilitated a new paradigm of topological deep learning. This approach has demonstrated superior efficacy across a range of biomolecular tasks, including protein-ligand binding affinity prediction \cite{meng2021persistent}, protein-protein interaction analysis \cite{wee2022persistent, liu2022hom}, and the early prediction of dominant SARS-CoV-2 variants, such as Omicron BA.4 and BA.5, two months
prior to their announcement by the World Health Organization (WHO) \cite{chen2022persistent}.

Despite their established utility, neither standard persistent homology nor the PL is able to capture heterogeneous information within complex datasets. To address this limitation, element-specific persistent homology was originally developed to explicitly distinguish between different atomic types \cite{cang2018representability}, an approach that subsequently inspired a wave of novel TDA methodologies \cite{ameneyro2024quantum, grbic2022aspects, liu2022biomolecular, liu2022hom}. Building on this foundation, Persistent Sheaf Laplacian (PSL) \cite{wei2025persistent} has recently emerged as a more mathematically elegant theory capable of embedding heterogeneous information, such as geometry and partial charges, directly into topological analysis via the theory of cellular sheaves \cite{hansen2019toward, curry2014sheaves}. Distinct from standard Laplacians, this sheaf-theoretic framework allows the assignment of specific vector spaces, representing chemical properties, to the open sets of protein topology \cite{wei2025persistent}. By associating specific weights with atoms (nodes), PSL encodes local topological and geometric information directly into its harmonic and non-harmonic spectra, thereby enabling the precise capture of intrinsic mutation-induced physical and chemical interactions. A review of persistent topological Laplacians is available in \cite{wei2025persistent2}. TDA approaches beyond persistent homology is surveyed in \cite{su2025topological}. Recent advances in TDL are reviewed in \cite{papamarkou2024position}. For application perspectives, the reader is further referred to a comprehensive review of TDA and TDL in molecular sciences \cite{wee2022persistent}.

In this work, we introduce SheafLapNet, a novel deep learning framework grounded in PSL theory, designed to predict mutation-induced changes in both protein stability and solubility. The SheafLapNet architecture synergizes three distinct categories of features to construct a comprehensive molecular representation: (1) multiscale topological and geometric features derived from Persistent Sheaf Laplacians, (2) auxiliary physicochemical descriptors capturing local atomic interactions, and (3) evolutionary sequence embeddings extracted from the pre-trained ESM-2 protein transformer \cite{lin2023evolutionary}. These diverse feature sets are integrated into a unified neural network framework that can predict both stability and solubility alterations. We rigorously evaluated SheafLapNet against established benchmarks, employing the S2648 and S350 datasets for stability prediction and the PON-Sol2 dataset for solubility classification. Across these comprehensive assessments, SheafLapNet consistently outperforms existing state-of-the-art models. These results demonstrate a predictive accuracy improvement, validating the efficacy of incorporating sheaf-theoretical invariants for mutation impact prediction.

\section{Results}
\subsection{Overview of SheafLapNet}

\begin{figure}[t] 
\centering  \includegraphics[width=1.0\linewidth, height=10cm]{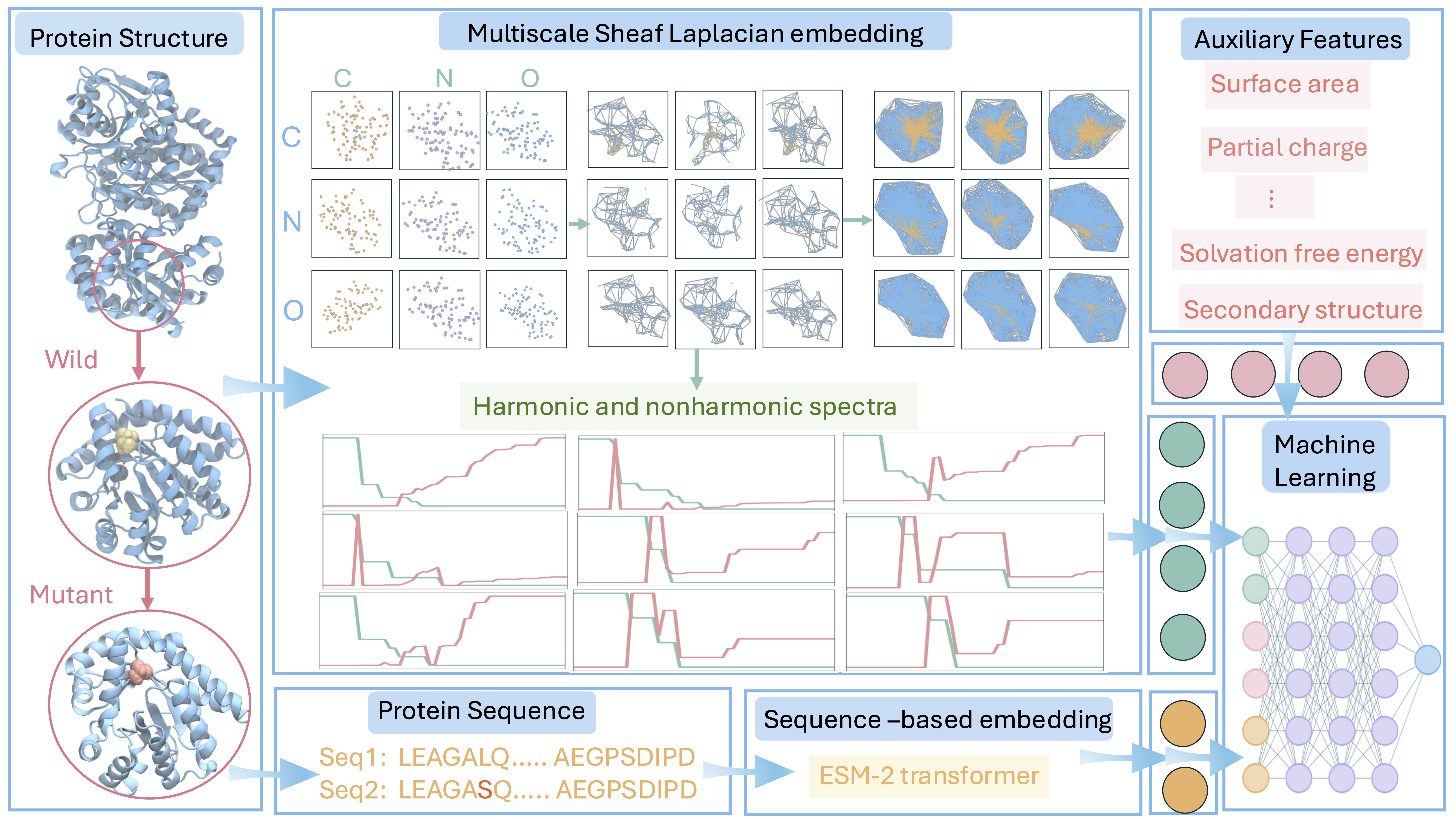} \caption{Illustration of the Persistent Sheaf Laplacian (PSL) neural network (SheafLapNet) workflow. The framework predicts mutation-induced stability and solubility changes by integrating multi-scale protein representations. For each input protein structure, the feature generation pipeline extracts three distinct components: (1) sequence-based embeddings derived from pretrained protein Transformer models, (2) topological features computed via the PSL framework, and (3) auxiliary physicochemical features. These three sets of features are concatenated to form the input of the neural network for the prediction task.} \label{fig:workflow} 
\end{figure}
Figure \ref{fig:workflow} outlines the workflow of SheafLapNet. As a standard machine learning model, SheafLapNet extracts features from protein structure and uses a neural network to predict mutation-induced stability and solubility changes upon mutation. The workflow begins with 3D protein structures from datasets, with corresponding mutant structures generated using the Jackal software \cite{xiang2002jackal}. The feature generation process consists of three components: sequence features from pretrained protein Transformer, topological features from persistent Sheaf Laplacians, and auxiliary physicochemical features. For topological features, atom subsets around the mutational site are extracted from both wild-type and mutant proteins to form element-specific subcomplexes. These subcomplexes are utilized to compute the harmonic and nonharmonic spectra of sheaf Laplacians under a structural filtration, creating a Sheaf Laplacian embedding that characterizes atom-atom interactions across multiple scales. For sequence features, the FASTA sequences of the wild-type
and mutant proteins are extracted from the complex and input into the pretrained Transformer models. The derived latent space embeddings are used as the sequence features. For physicochemical features, we consider the atom-level properties such as partial
charge, electrostatic solvation free energy, and Coulomb interactions, as well as residue-level
properties such as mutation site neighbor amino acid composition, pKa shifts, and additional
physicochemical properties. These three types of feature embeddings are concatenated to
form the input feature vector for the machine learning algorithm (neural network) to predict the mutation-induced stability changes and solubility changes.

\subsection{Prediction of Mutation-Induced Protein Stability Changes}
\begin{figure}[t] 
\centering  \includegraphics[width=1\linewidth, height=12cm]{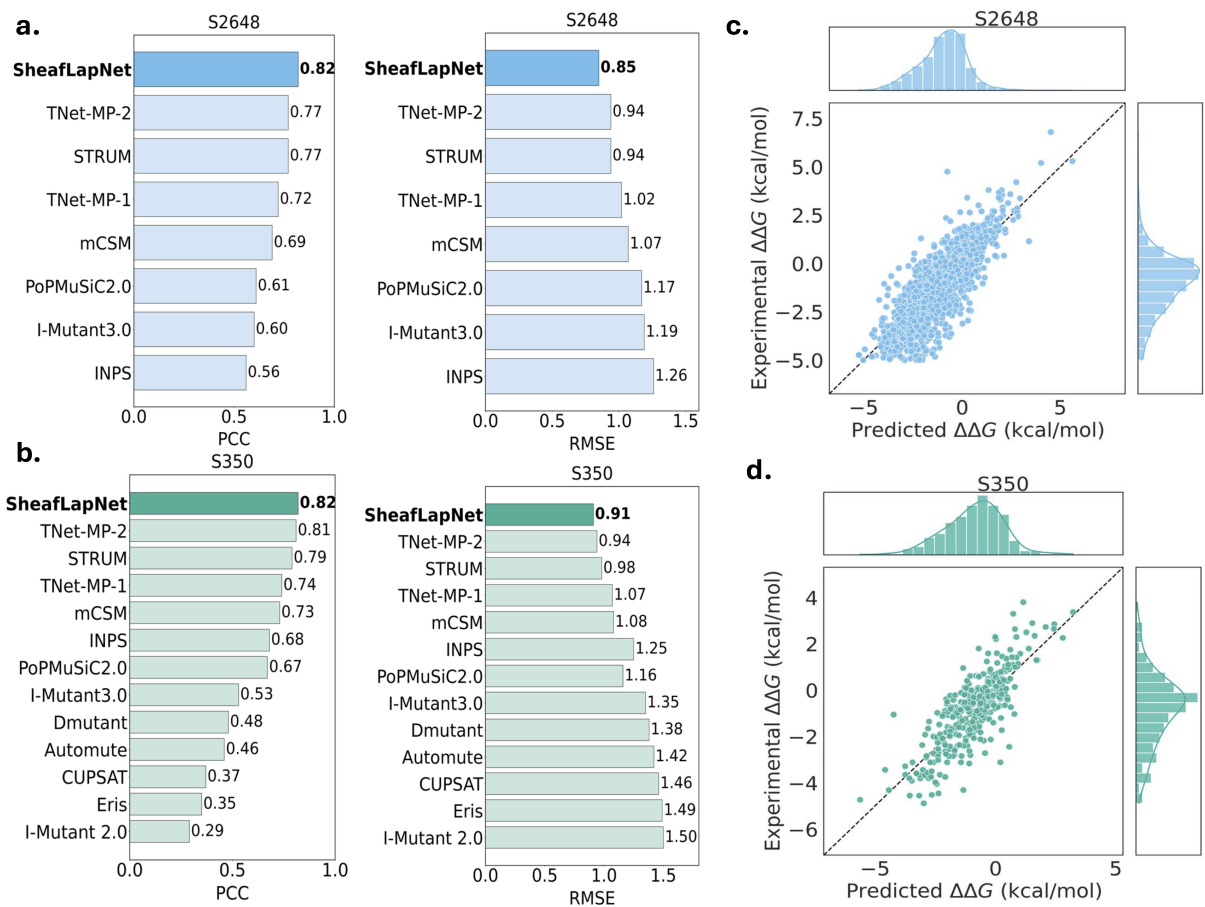} \caption{Illustration of model performance in protein stability changes upon mutation. (a). 5-fold cross-validation performance of SheafLapNet for S2648 dataset compared to existing state-of-the-art models \cite{cang2017topologynet, quan2016strum, worth2011sdm}.
(b). Blind test
performance of SheafLapNet for S350 dataset compared to existing state-of-the-art models \cite{cang2017topologynet, quan2016strum, worth2011sdm, capriotti2005mutant2}. (c). Comparison of experimental protein stability changes with predicted ones from SheafLapNet for
S2648 dataset. (d). Comparison of experimental protein stability changes with predicted ones from SheafLapNet for S350 dataset. } \label{fig:stability} 
\end{figure}

\begin{figure}[!htbp] 
\centering  \includegraphics[width=1\linewidth]{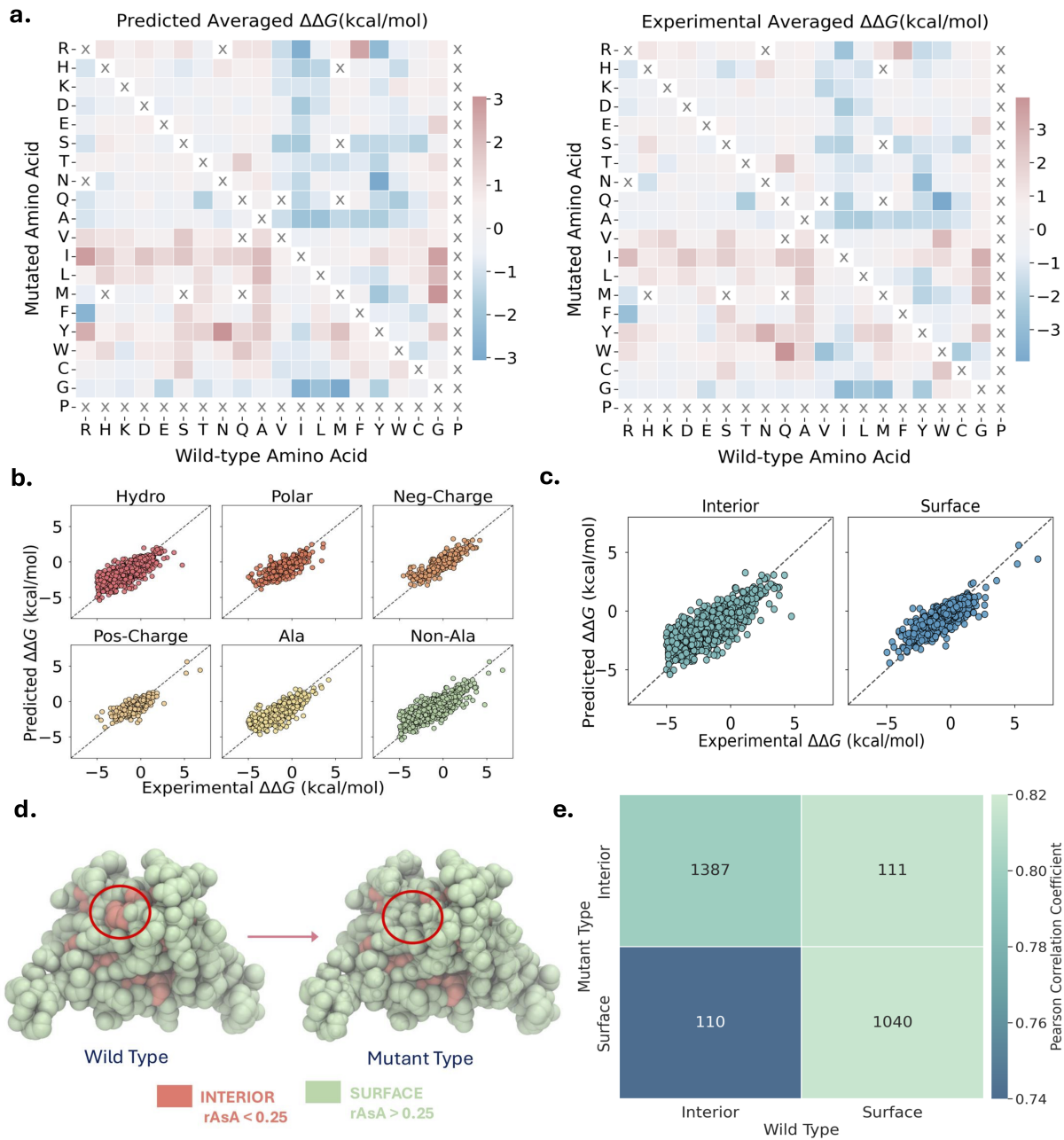}
 \caption{Illustration of model performance in protein stability changes upon mutation in the S2648 dataset. (a). Residue-residue matrix comparing the averaged experimental and
predicted mutation-induced stability changes $\Delta\Delta G$ for the whole dataset of 2648 mutations, where X indicates no mutation samples. (b). Model performance
across different mutation types based on residue physicochemical properties. (c): Model performance across
different structural regions defined by the relative accessible solvent area (rASA). (d): Structural shift on protein (PDB ID: 1ARR) of the Q39G mutation from interior to surface region. (e): Model performance stratified by four mutation region types, with cell annotations indicating the sample count per category.} \label{fig:heatmap} 
\end{figure}

Mutation-induced perturbations in protein stability are a fundamental mechanism underlying numerous genetic diseases. These stability alterations, quantified as $\Delta\Delta G$, frequently compromise protein function and structural integrity. To rigorously assess the capacity of our framework to capture these effects, we utilized the S2648 benchmark dataset, which comprises 2,648 single-point mutations across 131 protein structures annotated with experimentally determined stability changes. Our evaluation protocol consisted of two distinct phases designed to ensure robust validation. First, we performed a fivefold cross validation on the complete S2648 dataset to establish baseline generalizability. This was followed by a targeted blind test using the S350 dataset, a high quality and curated subset of S2648 comprising 350 mutations across 67 proteins. This subset serves as a standard benchmark for independent model evaluation. To facilitate direct comparison with established methods in the literature \cite{cang2017topologynet, quan2016strum, worth2011sdm}, predictive performance was quantified using the Pearson correlation coefficient (PCC) and root mean squared error (RMSE) as defined in Section S1.1 of the Supporting Information.

The predictive efficacy of the Sheaf Laplacian framework is substantiated by the superior performance of SheafLapNet across both validation protocols. In the comprehensive S2648 cross-validation analysis, we benchmarked our model against existing state-of-the-art methods as shown in Fig.~\ref{fig:stability}(a). SheafLapNet significantly outperforms established baselines, including STRUM \cite{quan2016strum} and the leading topological convolutional neural network, TNet-MP-2 \cite{cang2017topologynet}. While these methods achieve a PCC of approximately 0.77 and an RMSE of 0.94 kcal/mol, our model demonstrates superior predictive capability, attaining a PCC of 0.82 and an RMSE of 0.85 kcal/mol. This corresponds to a 6.49\% improvement in Pearson correlation and a 9.6\% reduction in prediction error relative to TNet-MP-2. As illustrated in Fig.~\ref{fig:stability}(c), the scatter plot demonstrates a tight concordance between predicted and experimentally measured stability changes, validating the efficacy of PSL embeddings to mutation-induced structural shifts. To further evaluate the robustness of our model, we conducted a blind test on the S350 benchmark. In this rigorous assessment, SheafLapNet maintained high predictive accuracy, achieving an average PCC of 0.82 and an RMSE of 0.91 kcal/mol. As depicted in the comparative analysis in Fig.~\ref{fig:stability}(b), our model outperforms TNet-MP-2, which reports a PCC of 0.81 and an RMSE of 0.94 kcal/mol. A direct comparison between the experimental and predicted $\Delta\Delta G$ values for the test set is presented in Fig.~\ref{fig:stability}(d). Beyond predictive accuracy, the PSL framework demonstrates superior robustness regarding data coverage. For instance, I-Mutant 3.0 failed to evaluate the complete dataset, covering only 2,636 of 2,648 samples in the S2648 benchmark and 338 of 350 in S350. In contrast, SheafLapNet successfully generated predictions for every mutation sample in both benchmarks. This complete coverage underscores the scalability of Sheaf Laplacian embeddings and highlights their computational reliability.

Using the predictions obtained from the 5-fold cross-validation on the S2648 dataset, we evaluate the concordance between average experimental and predicted stability changes across diverse mutation types. As shown in the Fig.~\ref{fig:heatmap}(a), a residue-residue matrix is established where the x-axis corresponds to the wild-type residues and the y-axis denotes the mutated residues. Each cell at the intersection of a specific mutant row and wild-type column represents the averaged $\Delta\Delta G$ value for that specific substitution pair, and the corresponding number of mutation samples is provided in Fig. S1 for each substitution. To capture the complete stability landscape, we explicitly incorporated reverse mutations by assigning opposite stability change values, resulting in an antisymmetric matrix structure. The predicted matrix closely mirrors the patterns observed in the experimental benchmark, indicating that the model successfully captures mutation-specific stability trends. While the overall landscape of perturbations is preserved, the variance in the predicted values appears attenuated compared to the experimental data. A prominent feature emerging from this analysis is that mutations to Isoleucine consistently yield positive $\Delta\Delta G$ values across the majority of wild-type residues. A possible explanation is that Isoleucine is often unable to compensate for the loss of favorable hydrophilic interactions, such as hydrogen bonds and electrostatic contacts, that are originally contributed by polar or charged wild-type residues. This leads to a lack of chemical compatibility with the surrounding aqueous solvent and consequently results in positive $\Delta\Delta G$ values.

To further evaluate the sensitivity of the model to specific physicochemical properties, mutation residues were stratified into six distinct categories: hydrophobic, polar, negatively charged, positively charged, alanine, and non-alanine. The correspondence between predicted and experimental stability changes for each category is presented in Fig.~\ref{fig:heatmap}(b). Quantitative analysis reveals consistent predictive performance across all groups. Specifically, the PCC (RMSE) values are 0.802 (0.890 kcal/mol), 0.750 (0.816 kcal/mol), 0.853 (0.768 kcal/mol), 0.798 (0.754 kcal/mol), 0.815 (0.859 kcal/mol), and 0.820 (0.832 kcal/mol) for hydrophobic, polar, negatively charged, positively charged, alanine, and non-alanine mutations, respectively. Notably, robust results were obtained across all mutation types, indicating that the model generalizes effectively across diverse amino acid substitutions without exhibiting significant bias. The consistently high correlation observed across these varying chemical groups suggests that the underlying feature representation effectively captures the distinct stabilizing and destabilizing forces associated with each specific residue type.

Mutation sites are stratified into two distinct structural environments: the protein interior and the solvent-exposed surface. This classification is determined by the relative accessible solvent area (rAsA), utilizing a threshold of 0.25. Residues satisfying $\text{rAsA} < 0.25$ are categorized as interior, while those with $\text{rAsA} \geq 0.25$ are classified as surface. Biologically, mutations occurring within the protein interior generally exhibit a larger magnitude of destabilization compared to surface mutations. In our dataset, the average experimental $\Delta\Delta G$ for interior residues is $-1.392$~kcal/mol, compared to $-0.522$~kcal/mol for surface residues. This trend confirms the greater destabilizing effect in buried regions, attributable to the disruption of the tightly packed hydrophobic core essential for protein stability. The predictive performance of the model across these structural regions is illustrated in Fig.~\ref{fig:heatmap}(c). The analysis reveals strong correlations in both environments: the interior region yields a PCC of 0.796 and an RMSE of 0.964~kcal/mol, while the surface region exhibits slightly superior performance with a PCC of 0.819 and an RMSE of 0.614~kcal/mol. This lower error on the surface likely reflects the lower variance in stability changes typical of solvent-exposed sites. To visualize these structural definitions, Fig.~\ref{fig:heatmap}(d) depicts the structural shift in the protein 1ARR induced by the Q39G mutation. In this instance, residue 39 transitions from a Glutamine in the interior region of the wild type to a Glycine in the surface region of the mutant. By categorizing residues as either interior or surface for both wild-type and mutant structures, we can examine the influence of continuous amino acid exposure on stability changes post-mutation. Figure~\ref{fig:heatmap}(e) displays the performance metrics across four specific mutation trajectories: Interior-to-Interior, Interior-to-Surface, Surface-to-Interior, and Surface-to-Surface. The model achieves robust performance across the majority of categories, with PCC values of 0.80 for Interior-to-Interior, 0.81 for Surface-to-Interior, and 0.82 for Surface-to-Surface. The slightly lower correlation observed for the Interior-to-Surface category with PCC value 0.74 is likely attributable to the significantly smaller sample size available for this specific mutation type.

\subsection{Classification of Mutation-Induced Protein Solubility Changes}
\begin{figure}[!htbp] 
\centering  \includegraphics[width=1\linewidth]{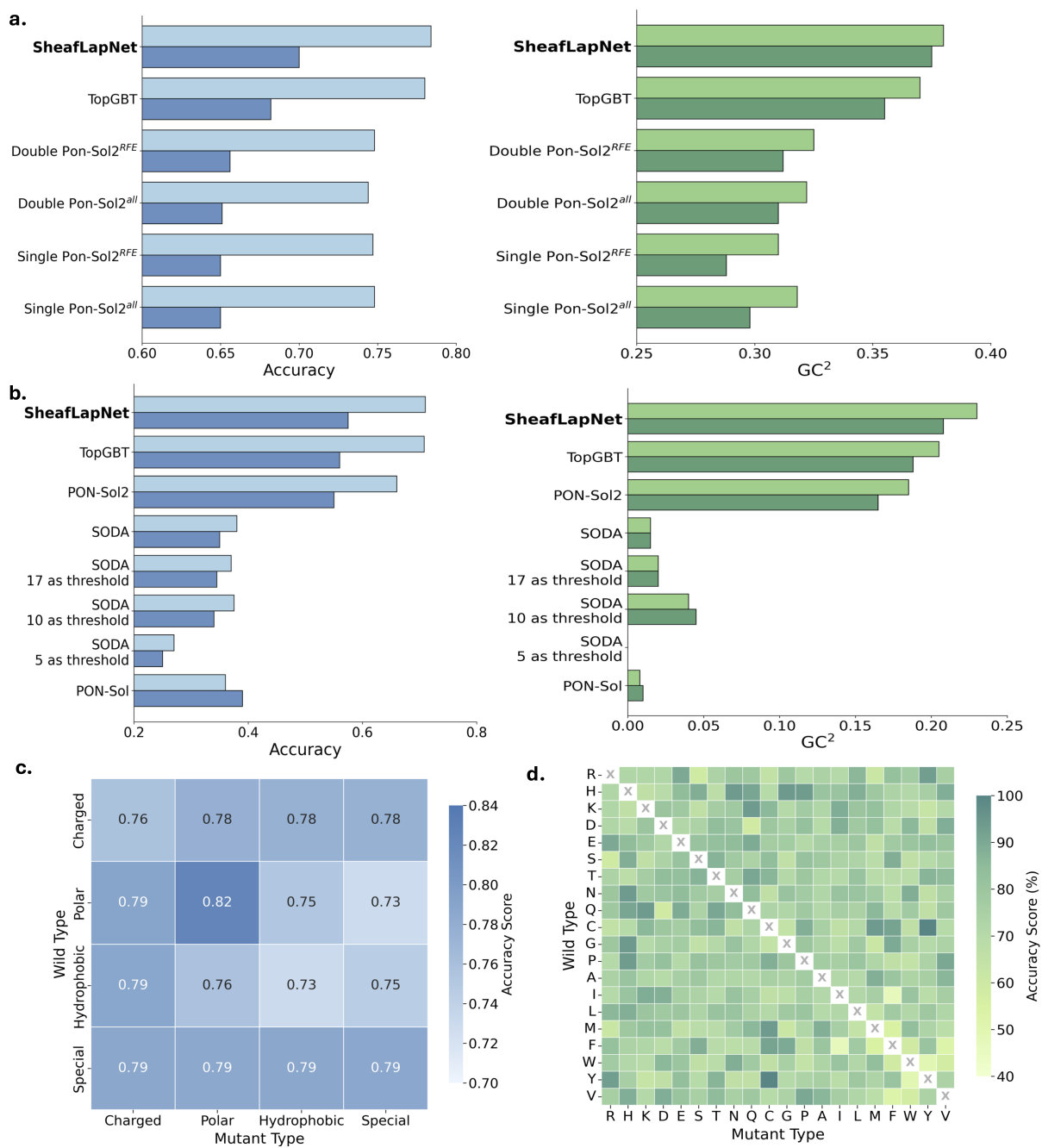} \caption{Illustration of model performance in classifying mutation-induced protein solubility changes.
Light blue denotes Accuracy (CPR), dark blue denotes Normalized Accuracy, light green denotes Generalized Squared Correlation (GC$^2$), and dark green denotes Normalized GC$^2$.
(a) Performance from 10-fold cross-validation on the PON-Sol2 dataset, comparing SheafLapNet against existing models \cite{yang2021pon, wee2024integration}. For existing PON-Sol2 models \cite{yang2021pon}, RFE refers to recursive feature elimination, and all refers to the use of all features.
(b) Independent blind test results on the PON-Sol2 dataset, comparing SheafLapNet against state-of-the-art models \cite{yang2021pon, yang2016pon, wee2024integration}.
(c) Model accuracy from 10-fold cross-validation stratified by physicochemical mutation groups.
(d) Model accuracy from 10-fold cross-validation analyzed by specific amino acid substitution types, where X denotes samples with no mutation.} \label{fig:mutsol} 
\end{figure}


Changes in protein solubility arising from genetic mutations are critical determinants of protein stability and function, often leading to aggregation-related diseases or altered biological activity. To model these effects, we utilized the dataset originally curated for PON-Sol2 \cite{yang2021pon}. This dataset comprises 6,328 mutation samples derived from 77 distinct proteins. Mutations are categorized into three classes based on their impact on solubility: decrease, increase, and neutral (no change). The distribution and classification of these mutation samples are illustrated in Fig. S2. Specifically, the dataset contains 3,136 samples demonstrating a decrease, 1,026 samples indicating an increase, and 2,166 samples showing no change. This distribution reveals a notable class imbalance with a ratio of approximately 1 : 0.69 : 0.33 (Decrease : Neutral : Increase), indicating a prevalence of mutations that reduce protein solubility.

To rigorously assess the predictive capabilities of our framework, we implemented a dual validation strategy comprising a random 10-fold cross-validation and an independent blind test classification. In the cross-validation phase, SheafLapNet was benchmarked against TopGBT \cite{wee2024integration}, which is a model grounded in persistent homology, and the existing PON-Sol2 predictive models\cite{yang2021pon}, which employ feature selection techniques such as recursive feature elimination(RFE). To mitigate the effects of the inherent class imbalance, normalized accuracy scores were utilized as the primary metric for comprehensive performance assessment. Given that we are dealing with a K-class problem with three distinct solubility classes, we rely on the Correct Prediction Ratio (CPR) and Generalized Squared Correlation ($\text{GC}^2$) as defined in Section S1.2 of the Supporting Information to provide a holistic assessment. Specifically, CPR measures the overall accuracy of the model while GC$^2$ quantifies the correlation coefficient of the classification, ranging from 0 to 1. Larger values for these metrics denote better performance.

The proposed model, SheafLapNet, achieved normalized CPR and GC$^2$ scores of 0.700 and 0.375, respectively, in the 10-fold cross-validation setting. These results surpass the existing PON-Sol2 models by up to 6.71\% and 20.19\%, and the persistent homology-based TopGBT model by up to 2.64\% and 5.63\%. Given the multi-class nature of the task, these metrics provide a comprehensive assessment of predictive performance, as illustrated in Fig.~\ref{fig:mutsol}(a). The robustness of the Sheaf Laplacian framework was further confirmed through independent blind testing. In this evaluation, SheafLapNet maintained consistent predictive power, achieving a normalized CPR of 0.570 and a normalized GC$^2$ of 0.201. These scores exceed those of the PON-Sol2 models by up to 3.64\% and 21.82\%, and TopGBT by up to 1.79\% and 6.91\%. Benchmark comparisons based on CPR and GC$^2$ scores are presented in Fig.~\ref{fig:mutsol}(b).


Switching focus to mutation types, our model's capability in classifying solubility changes also merits exploration across the 20 distinct amino acid types in the dataset. In addition to this, we subgroup amino acids as charged polar, hydrophobic, or special cases. Fig.~\ref{fig:mutsol}(c) displays accuracy scores for each mutation group pair, while Fig.~\ref{fig:mutsol} (d) shows scores for each amino acid pair. In this representation, the y-axis labels the residue type of the original protein, whereas the x-axis labels the residue type of the mutant. For a reverse mutation, the labels are taken with reverse solubility change unless the change is zero. The analysis reveals a clear disparity in predictive performance across mutation types. Notably, mutations transitioning from polar residues to polar groups achieve the highest classification accuracy. This is likely due to the direct correlation between surface polarity and solubility, resulting in predictable solvation energy changes that are effectively captured by the model's topological features. Conversely, transitions from polar to special residues or within the hydrophobic group exhibit comparatively lower performance. This reduced accuracy may reflect the challenge of modeling complex structural constraints, such as backbone rigidity introduced by special residues, and the subtle packing interactions characteristic of the hydrophobic core.

\section{Materials and Methods}

This section details the three distinct feature sets extracted from protein structures for our analysis: topological features derived from Persistent Sheaf Laplacians, physicochemical auxiliary descriptors, and evolutionary sequence embeddings. A comprehensive summary of the software packages utilized in this study is provided in Section S4 of the Supplementary Information, while a detailed list of the auxiliary descriptors can be found in Section S3.

\subsection{Persistent Sheaf Laplacian}

Before introducing the Persistent Sheaf Laplacian, we first review the standard Persistent Laplacian. The construction begins by considering two simplicial complexes, $K$ and $L$, such that $K \subseteq L$. We denote the simplicial chain complexes of $K$ and $L$ with real coefficients as $C^K$ and $C^L$, respectively. Since a chain group $C_q$ in a simplicial chain complex is formally generated by simplices, it naturally forms a finite-dimensional inner product space, ensuring that the adjoint of the boundary map $\partial_q$ is well-defined. 

We define a subspace $C^{L,K}_{q+1}$ of $C^L_{q+1}$ as the set $\{c \in C^L_{q+1} \mid \partial^L_{q+1}(c) \in C^K_q\}$. Let $\partial^{L,K}_{q+1}$ denote the restriction of the boundary map $\partial^L_{q+1}$ to this subspace $C^{L,K}_{q+1}$. The $q$-th Persistent Laplacian $\Delta^{L,K}_q$ is formally defined by the operator:
\begin{equation}
    \Delta^{L,K}_q = \partial^{L,K}_{q+1}(\partial^{L,K}_{q+1})^* + (\partial^K_q)^*\partial^K_q.
\end{equation}

To extend this framework to the Persistent Sheaf Laplacian, we must first define a cellular sheaf. A cellular sheaf $\mathscr{S}$ consists of a simplicial complex $X$ (viewed as a cell complex) equipped with an assignment to each cell $\sigma$ of $X$ a finite-dimensional vector space $\mathscr{S}(\sigma)$, referred to as the stalk of $\mathscr{S}$ over $\sigma$. Additionally, for each face relation $\sigma \preceq \tau$ (where $\sigma \subset \bar{\tau}$), there exists a linear morphism of vector spaces denoted by $\mathscr{S}_{\sigma \preceq \tau}$, known as the restriction map. These maps satisfy the composition rule:
\begin{equation}
    \mathscr{S}_{\rho \preceq \tau} = \mathscr{S}_{\sigma \preceq \tau}\mathscr{S}_{\rho \preceq \sigma}
\end{equation}
for any $\rho \preceq \sigma \preceq \tau$, with $\mathscr{S}_{\sigma \preceq \sigma}$ being the identity map of $\mathscr{S}(\sigma)$.

Analogous to a simplicial complex, a cellular sheaf generates a sheaf cochain complex. The $q$-th sheaf cochain group $C^q_{\mathscr{S}}$ is defined as the direct sum of stalks over $q$-dimensional cells. By globally orienting the simplicial complex $X$ to obtain a signed incidence relation $[\sigma : \tau]$ for each $\sigma \preceq \tau$, the coboundary map $d^q : C^q_{\mathscr{S}} \to C^{q+1}_{\mathscr{S}}$ is defined by:
\begin{equation}
    d^q|_{\mathscr{S}(\sigma)} = \sum_{\sigma \preceq \tau} [\sigma : \tau] \mathscr{S}_{\sigma \preceq \tau}.
\end{equation}

The Persistent Sheaf Laplacian is constructed by considering two cellular sheaves, $\mathscr{S}$ on $K$ and $\mathscr{T}$ on $L$, such that $K \subseteq L$ and the stalks and restriction maps of $K$ are identical to those of $L$. The relationship between the cochain complexes of these sheaves and the persistent spectral structures is illustrated in the following commutative diagram:

\begin{equation}
\begin{tikzcd}[row sep=3.5em, column sep=4.5em]
    C^{q-1}_{\mathscr{S}} 
    \arrow[r, shift left=1ex, "d^{q-1}_{\mathscr{S}}"] 
    & C^q_{\mathscr{S}} 
    \arrow[l, shift left=1ex, "(d^{q-1}_{\mathscr{S}})^*"] 
    \arrow[rr, "d^q_{\mathscr{S}}"] 
    \arrow[dr, shift left=0.5ex, "d^q_{\mathscr{S},\mathscr{T}}"] 
    & 
    & C^{q+1}_{\mathscr{S}} 
    \\
    & 
    & C^{q+1}_{\mathscr{S},\mathscr{T}} 
    \arrow[ul, shift left=0.5ex, "(d^q_{\mathscr{S},\mathscr{T}})^*"] 
    \arrow[dr, dashed, hook] 
    & 
    \\
    C^{q-1}_{\mathscr{T}} 
    \arrow[uu, "\pi"] 
    \arrow[r, "d^{q-1}_{\mathscr{T}}"] 
    & C^q_{\mathscr{T}} 
    \arrow[uu, "\pi"] 
    \arrow[rr, "d^q_{\mathscr{T}}"] 
    & 
    & C^{q+1}_{\mathscr{T}} 
    \arrow[uu, "\pi"]
\end{tikzcd}
\end{equation}

We define the subspace $C^{q+1}_{\mathscr{S},\mathscr{T}}$ as $ C^{q+1}_{\mathscr{S},\mathscr{T}} = \{c \in C^{q+1}_{\mathscr{T}} \mid (d^q_{\mathscr{T}})^*(c) \in C^q_{\mathscr{S}}\}.$
Let $d^q_{\mathscr{S},\mathscr{T}}$ denote the adjoint map of the restricted operator $(d^q_{\mathscr{T}})^*|_{C^{q+1}_{\mathscr{S},\mathscr{T}}}$. The $q$-th PSL, denoted as $\Delta^{\mathscr{S},\mathscr{T}}_q$, is defined as:
\begin{equation}
    \Delta^{\mathscr{S},\mathscr{T}}_q = (d^q_{\mathscr{S},\mathscr{T}})^* d^q_{\mathscr{S},\mathscr{T}} + d^{q-1}_{\mathscr{S}} (d^{q-1}_{\mathscr{S}})^*.
\end{equation}
The spectral properties of the PSL provide key topological invariants. Specifically, the dimension of the zero-eigenspace (kernel) of the operator corresponds to the $q$-th persistent sheaf Betti number, denoted $\beta_q^{\mathscr{S},\mathscr{T}}$ with $\beta_q^{\mathscr{S},\mathscr{T}} = \dim(\ker \Delta_q^{\mathscr{S},\mathscr{T}})$. 

This study utilizes cellular sheaves constructed on labeled simplicial complexes to model protein structures using atomic partial charges, following methodologies established in existing literature \cite{wei2025persistent}.
A general framework is first defined for constructing sheaves on a labeled simplicial complex $X$, where each vertex is associated with a quantity $q$. Let $F: X \to \mathbb{R}$ be a nowhere-zero function. The sheaf is constructed by assigning the stalk $\mathbb{R}$ to each simplex. For a face relation $[v_0, \dots, v_n] \leq [v_0, \dots, v_n, v_{n+1}, \dots, v_m]$ (where orientation is not relevant), the linear morphism 
$\mathscr{S}([v_0, \dots, v_n] \leq [v_0, \dots, v_n, v_{n+1}, \dots, v_m])$
is defined as the scalar multiplication by:
\begin{equation}
    \frac{F([v_0, \dots, v_n]) q_{n+1} \cdots q_m}{F([v_0, \dots, v_n, v_{n+1}, \dots, v_m])}.
\end{equation}

\begin{figure}[t] 
\centering  \includegraphics[width=1\linewidth]{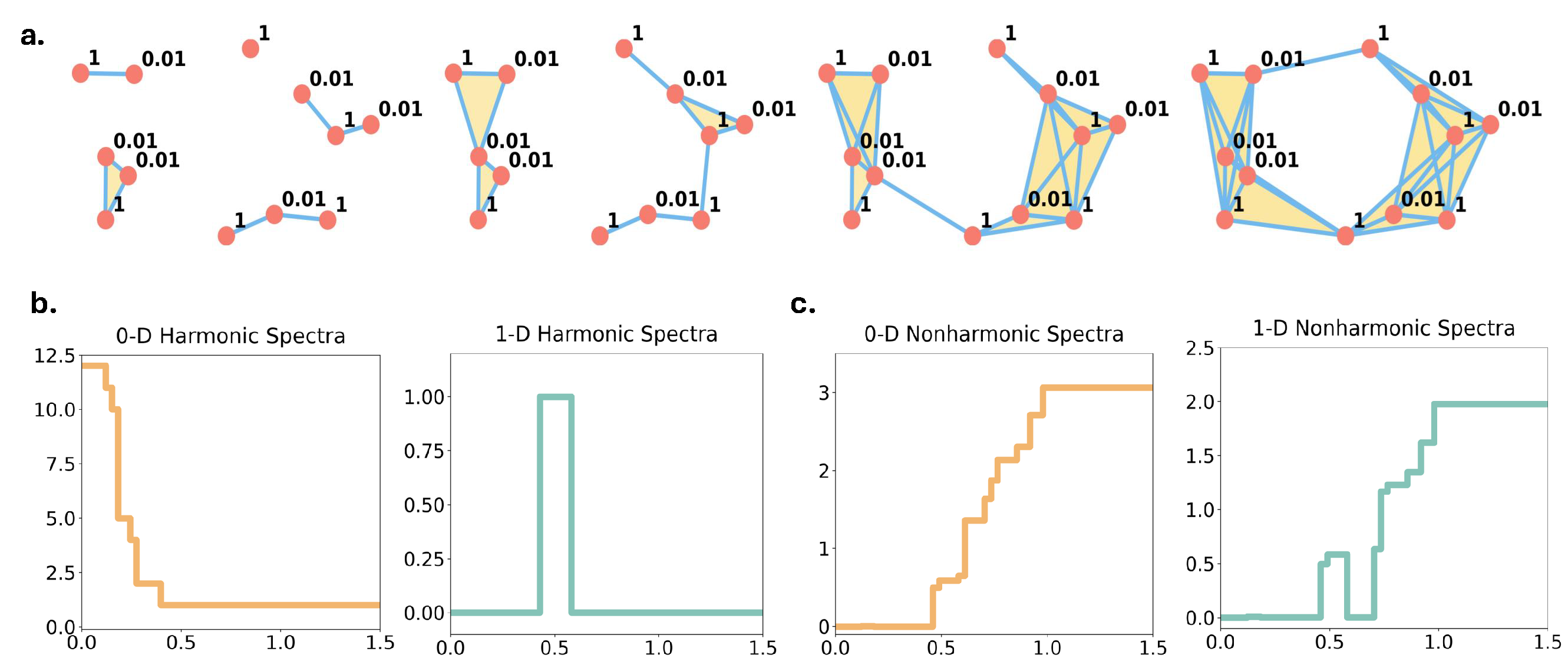}
\caption{ Illustration of persistent sheaf Laplacians. (a): a filtration process of the Rips complex from a point
cloud data. (b): the persistent multiplicity of zero eigenvalues from the sheaf Laplacian matrices in 0-D and 1-D dimensions.
(c): the minimal value of non-zero eigenvalues from the zero-dimensional and one-dimensional sheaf Laplacian matrices. } \label{fig:psl} 
\end{figure}

To adapt this framework for protein analysis, non-geometrical information is incorporated by employing atomic partial charges obtained from the PDB2PQR package \cite{jurrus2018improvements}. A Rips or Alpha filtration of graphs is constructed wherein vertices $v_i$ represent atoms and edges $e_{ij}$ represent interactions between atoms $v_i$ and $v_j$. 
The cellular sheaf is defined such that each stalk is the real line $\mathbb{R}$. For the specific face relation $v_i \preceq e_{ij}$ (where a vertex $v_i$ is a face of the edge $e_{ij}$), the restriction morphism is explicitly defined as multiplication by $\frac{q_j}{r_{ij}}$,
where $q_j$ is the partial charge of the neighboring atom $v_j$, and $r_{ij}$ represents the Euclidean length of the edge $e_{ij}$. The harmonic spectra of the resulting PSL reveal topological invariants, while the nonharmonic spectra represent geometric information of the data. Consequently, these spectra are utilized as the input features for the SheafLapNet model.

To illustrate the Persistent Sheaf Laplacian framework, we present an example using a point cloud dataset consisting of 12 points, as depicted in Figure~\ref{fig:psl}. The topological structure is modeled through a Vietoris-Rips (VR) complex filtration, visualized in Figure~\ref{fig:psl}(a). Here, heterogeneous information is explicitly mapped onto the domain, with nodes assigned distinct "charge" values of 1 and 0.01. As the filtration parameter increases, the connectivity of the point cloud evolves, generating a sequence of nested simplicial complexes that serves as the domain for our sheaf-theoretic analysis. Figure~\ref{fig:psl}(b) presents the corresponding harmonic spectral analysis derived from the persistent sheaf Laplacians, specifically plotting the persistent multiplicities of the zero eigenvalues. These multiplicities correspond to the Sheaf Betti numbers, which generalize standard topological invariants by incorporating heterogeneous information such as atomic partial charges. Complementing this harmonic analysis, Figure~\ref{fig:psl}(c) displays the non-harmonic spectral information, specifically the evolution of the minimum non-zero eigenvalues. These values quantify the connectivity strength and local geometric rigidity, offering a continuous measure of the intrinsic physical and chemical interactions encoded within the sheaf.

\subsection{Sheaf Laplacian feature generation for protein}
In the SheafLapNet architecture, protein structures are modeled as sets of simplicial complexes, from which PSLs are computed to extract topological features. To balance computational efficiency with the capability to capture essential physical interactions, the model employs an element-specific and site-specific atom subsetting strategy. Atoms within the three-dimensional protein structure are partitioned into mutation-site atoms, denoted as $A_m$, and mutation neighborhood atoms located within a cutoff radius $r$, denoted as $A_{mn}(r)$. This partitioning is applied to both wild-type and mutant structures. Leveraging the element-specific representation framework \cite{cang2018integration}, the model focuses on interactions involving carbon (C), nitrogen (N), and oxygen (O) atoms. By considering the intersections of these element types between the mutation site and its neighborhood, nine distinct pairwise atom combinations are constructed. These combinations encode specific interaction types. For instance, the pairing of carbon atoms in the mutation site ($A_C \cap A_m$) with those in the neighborhood ($A_C \cap A_{mn}(r)$) captures hydrophobic C-C interactions, whereas combinations involving carbon and oxygen ($A_C \cap A_m$ and $A_O \cap A_{mn}(r)$) characterize polar C-O interactions. These element-specific and site-specific sets subsequently underpin the multiscale sheaf Laplacian embeddings. 



To accurately capture the topological interactions within the defined atom sets, we construct VR complexes \cite{vietoris1927hoheren} and Alpha complexes \cite{edelsbrunner2011alpha}. For the VR complex construction, we employ a modified distance function, $d$, designed to specifically isolate interfacial interactions between distinct atom sets (e.g., between the mutation site $A_m$ and the neighborhood $A_{mn}(r)$). The modified metric is defined as:

\begin{equation}
    d(a_i, a_j) = 
    \begin{cases} 
        E_d(a_i, a_j) & \text{if } (a_i \in A_m \land a_j \in A_{mn}(r)) \lor (a_i \in A_{mn}(r) \land a_j \in A_m), \\
        \infty & \text{otherwise},
    \end{cases}
\end{equation}
where $E_d(a_i, a_j)$ denotes the Euclidean distance between atoms $a_i$ and $a_j$. By assigning an infinite distance to atom pairs belonging to the same set, this metric effectively filters out intra-set edges, thereby focusing the topological analysis exclusively on the interactions between the mutation site and its surrounding environment. Complementing this, Alpha complexes \cite{edelsbrunner2011alpha} are constructed using the standard Euclidean distance on the same atom sets to capture local geometric constraints.

\begin{figure}[t]
\centering
\includegraphics[width=0.8\linewidth, height=8cm]{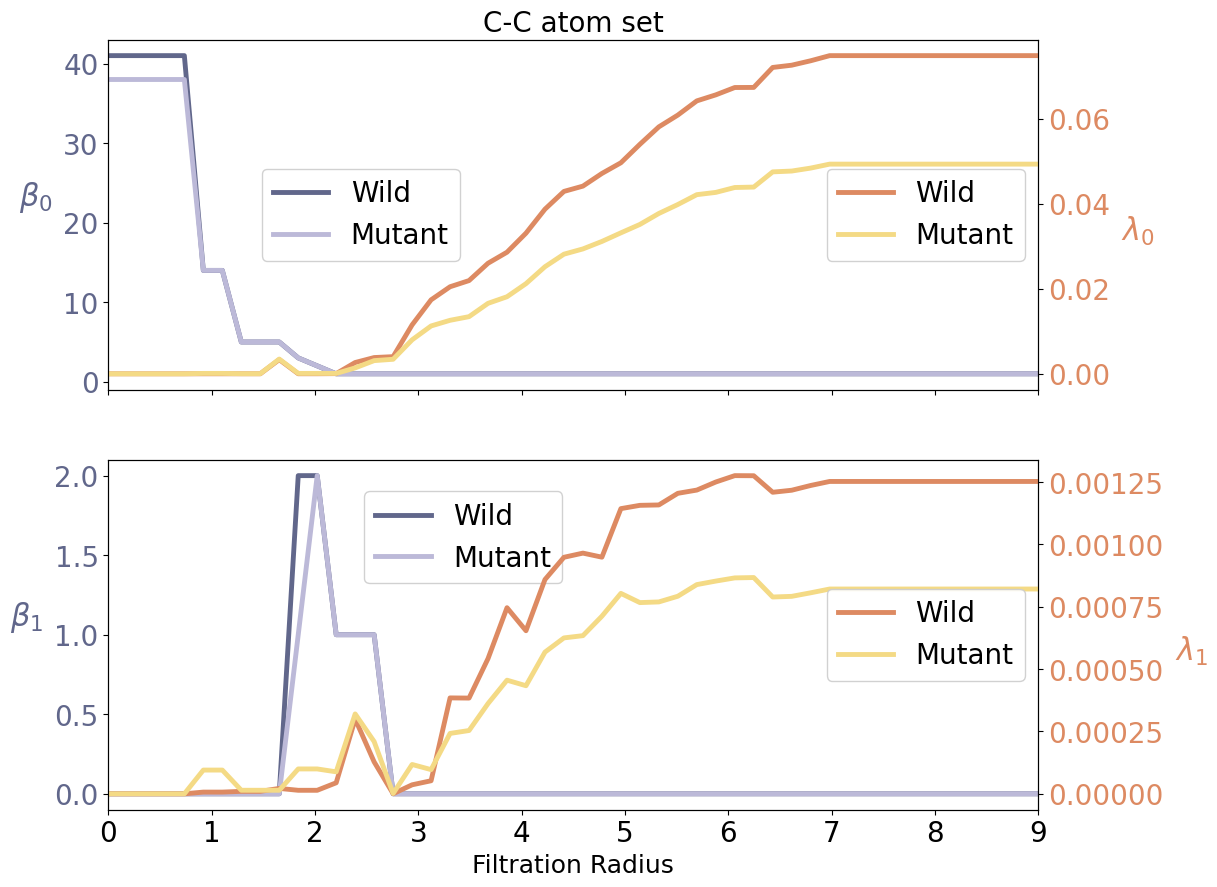}
\caption{An illustration of PSL descriptors for C-C interactions from VR complexes using modified distance filtration at residue 37 mutation site from L to S in protein 1A5E of S2648 dataset. Top panel: Illustration of zero-dimensional PSL descriptors. Bottom panel: Illustration of one-dimensional PSL descriptors. The left axis represents the Sheaf Betti numbers $\beta_0$ and $\beta_1$. The right axis represents the minimal nonzero eigenvalues $\lambda_0$ and $\lambda_1$ of the persistent sheaf Laplacian.}
\label{fig:barcode}
\end{figure}

\begin{figure}[t] 
\centering  \includegraphics[width=0.8\linewidth, height=8cm]{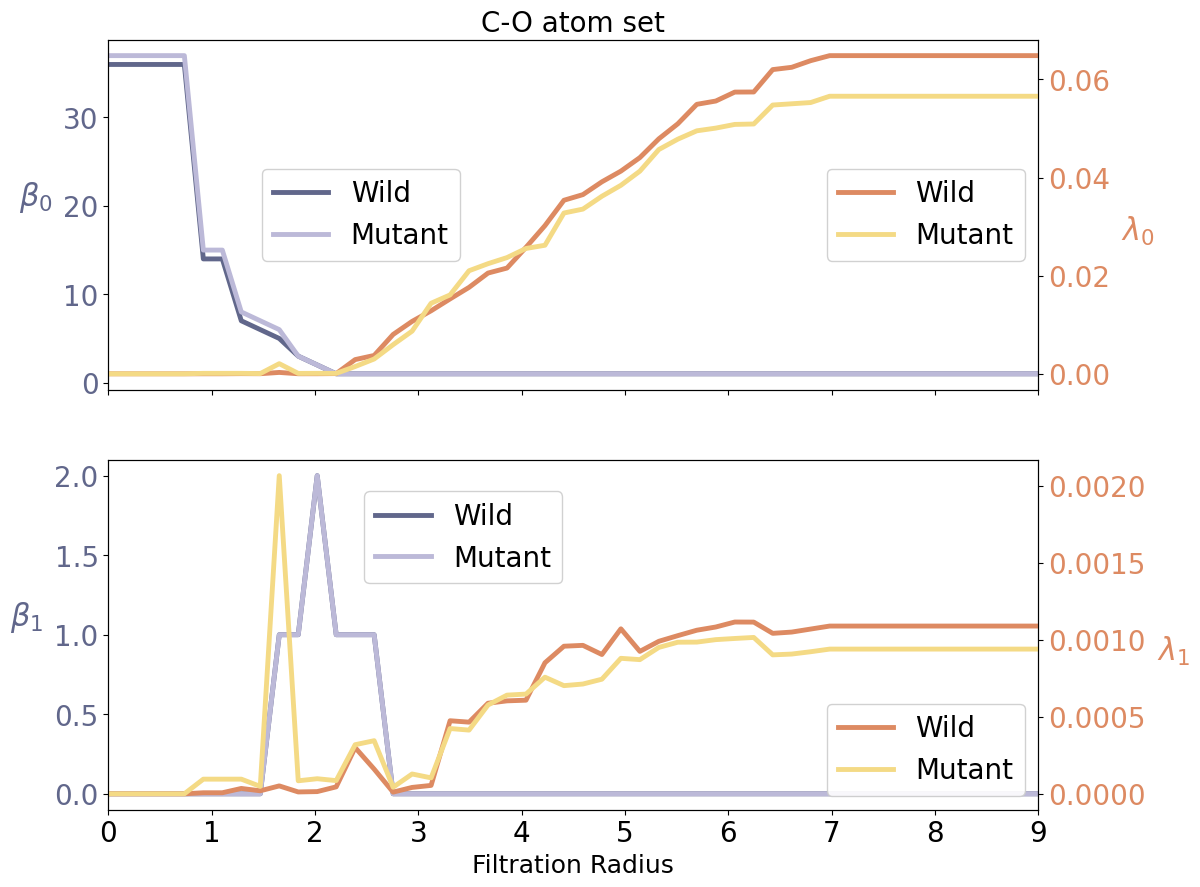} \caption{An illustration of PSL descriptors for C-O interactions from VR complexes using modified distance filtration at residue 37 mutation site from L to S in protein 1A5E of S2648 dataset. Top panel: Illustration of zero-dimensional PSL descriptors. Bottom panel: Illustration of one-dimensional PSL descriptors. The left axis represents the Sheaf Betti numbers $\beta_0$ and $\beta_1$. The right axis represents the minimal nonzero eigenvalues $\lambda_0$ and $\lambda_1$ of the persistent sheaf Laplacian.} \label{fig:barcode1} 
\end{figure}

Spectral features are extracted from the constructed complexes to serve as local PSL descriptors. In this study, a cutoff distance of $16\text{\AA}$ from the mutation site is employed to identify mutation neighborhood atoms. For the VR complex, the filtration is conducted over a range of $3\text{\AA}$ to $9\text{\AA}$ with a step size of $1\text{\AA}$. For each filtration step, spectral properties are derived from the Persistent Sheaf Laplacian. Regarding the harmonic components, the count of zero eigenvalues is recorded, generating a 7-dimensional feature vector for each atom set. For the non-harmonic components, eight statistical properties are extracted from the non-zero eigenvalue spectrum: maximum, minimum, mean, sum, standard deviation, variance, and the eigenvalue count. This procedure yields a 72-dimensional feature vector per atom set. Similarly, for the Alpha complex, the one-dimensional harmonic and non-harmonic components of the persistent sheaf Laplacians are analyzed using the same statistical extraction method. The final topological feature representation is constructed by combining the extracted features from the zero-dimensional VR models and the one-dimensional Alpha models. These feature vectors are computed for the wild-type structure, the mutant structure, and the difference between them. The concatenation of these vectors results in a comprehensive feature representation for a single protein, characterized by a 3402-dimensional vector. To illustrate the discriminative capacity of these features, Fig.~\ref{fig:barcode} and Fig.~\ref{fig:barcode1} depict the zero-dimensional and one-dimensional PSL descriptors for both wild-type and mutant structures of protein 1A5E from the S2648 dataset. Fig.~\ref{fig:barcode} utilizes atom sets $A_C \cap A_m$ and $A_C \cap A_{mn}(r)$ to generate VR complexes with a modified distance $d$-based filtration, effectively revealing hydrophobic C-C interactions. Similarly, Fig.~\ref{fig:barcode1} utilizes atom sets $A_C \cap A_m$ and $A_O \cap A_{mn}(r)$ to generate VR complexes with $d$-based filtration, revealing the polar C-O interactions.

\subsection {Auxiliary Features}

In addition to topological features, our model integrates a comprehensive set of physicochemical descriptors at both the atom and residue levels to characterize the mutation environment. At the atom level, we compute solvent-excluded surface areas, partial atomic charges, Coulombic and van der Waals interaction energies, and electrostatic solvation free energies \cite{chen2011mibpb, rocchia2001extending, jurrus2018improvements}. Here, electrostatic solvation free energies were computed using the MIBPB solver for the Poisson-Boltzmann model \cite{chen2011mibpb}.  At the residue level, the feature set comprises the amino acid composition and physical properties of the local neighborhood, mutation-induced pKa shifts, evolutionary conservation scores derived from Position-Specific Scoring Matrices (PSSM), and predicted secondary structure parameters. Detailed definitions and calculation protocols are provided in Section S3 of the Supporting Information.

\subsection{Sequence Features}

Protein large language models (pLLMs), such as ESM (Evolutionary Scale Modeling) \cite{rives2021biological} and ProtTrans \cite{elnaggar2021prottrans}, have demonstrated remarkable efficacy in capturing deep evolutionary patterns from hundreds of millions of unannotated sequences. To incorporate this rich evolutionary context into our algebraic framework, we employ the ESM-2 Transformer architecture \cite{lin2023evolutionary}. Trained via self-supervised masked language modeling (MLM) on a comprehensive protein sequence database, this architecture comprises 33 Transformer layers and approximately 650 million parameters. For a given input sequence, the model produces contextualized token embeddings of dimension 1,280. We extract these representations from the final layer and compute the mean across the sequence length, yielding a single fixed-length 1,280-dimensional vector. To explicitly capture the impact of genetic variations, we generate these embeddings for both the wild-type and mutant sequences; these are subsequently concatenated to form a final 2,560-dimensional evolutionary feature vector used in our predictive model.

\subsection{SheafLapNet model hyperparameters}

The SheafLapNet is a deep neural network designed to handle larger datasets by processing concatenated vectors of extracted PSL features, Transformer embeddings, and auxiliary vectors. 
The network architecture consists of six fully connected hidden layers, with each layer containing 15,000 neurons. To introduce non-linearity and enhance representational power, Rectified Linear Unit (ReLU) activation functions are applied after each hidden layer. Dropout regularization is incorporated throughout the network to prevent overfitting and ensure robust generalization. The final output layer is task-dependent, designed to generate a continuous scalar value for regression tasks in stability prediction or discrete class probabilities for classification tasks in solubility prediction. The model training was conducted over 200 epochs with a fixed learning rate of $0.001$. The batch size was optimized specifically for each task type: a batch size of 32 was employed for the regression-based protein stability datasets, while a batch size of 50 was utilized for the classification-based protein solubility datasets.

\section{Conclusion}

Accurately predicting the effects of genetic mutations on protein stability and solubility is fundamental to understanding the molecular basis of disease and advancing precision medicine. However, experimental characterization of these widely occurring variants is resource-intensive, and existing computational models often suffer from limited interpretability or treat stability and solubility as disjoint tasks. Therefore, developing unified, rigorous predictive frameworks capable of capturing intrinsic molecular interactions is a critical imperative. To address this challenge, we introduced SheafLapNet, a novel topological deep learning (TDL) framework that leverages the persistent sheaf Laplacian (PSL) to predict mutation-induced changes in both protein stability and solubility. By synergizing the heterogeneous geometric information encoded by PSL with auxiliary physicochemical descriptors and evolutionary sequence embeddings from the pre-trained ESM-2 transformer, our model achieves a mathematically grounded and robust representation of molecular alterations. This integrative approach captures the multiscale and mechanistic nuances of mutation impacts, resulting in a state-of-the-art predictive accuracy improvement across four rigorous benchmark datasets.

While this study establishes the efficacy of embedding heterogeneous atomic information into simplicial complexes via PSL, the rapidly evolving field of topological data analysis offers further opportunities for refinement. Future work could explore more distinct mathematical invariants, such as the persistent Path Laplacian \cite{wang2023persistent}, persistent hyperdigraph Laplacians \cite{chen2023persistent},  the persistent Dirac operator \cite{ ameneyro2024quantum, suwayyid2024persistent} and sheaf or copresheaf neural networks \cite{barbero2022sheaf,hajij2025copresheaf}. These advanced tools may provide deeper insights into the directed and higher-order interactions within protein structures that simplicial complexes alone may not fully capture, potentially unlocking new levels of predictive precision.

\section*{Data and Code Availability}
The data and code in this study can be found in \url{https://github.com/yren24/SheafLapNet/tree/main}

\section*{Acknowledgments}
This work was supported in part by NIH grant R35GM148196, National Science Foundation grant DMS2052983, and  Michigan State University Research Foundation.

\newpage 

\bibliographystyle{abbrv}

\bibliography{reference}

\end{document}